\documentclass{amsart}
\usepackage[cp1250]{inputenc}

\newtheorem{theorem}{Theorem}[section]

\newtheorem{example}[theorem]{Example}
\newtheorem{remark}[theorem]{Remark}
\newtheorem{conjecture}[theorem]{Conjecture}
\numberwithin{equation}{section}

\author{M. Primc}
\address[M. Primc]{Faculty of Science,  University of Zagreb,  Zagreb, Croatia}
\email{primc@math.hr}

\begin{document}

\title
 {New partition identities for odd $W$ odd}

\begin{abstract}
In this note we conjecture Rogers-Ramanujan type colored partition identities for an array ${\mathcal N}_w^{\text{odd}}$ with odd number of rows $w$ such that the first and the last row consist of even positive integers. In a strange way this is different from the partition identities for the array ${\mathcal N}_w$ with odd number of rows $w$ such that the first and the last row consist of odd positive integers---the partition identities conjectured by S. Capparelli, A. Meurman, A. Primc and the author and related to standard representations of the affine Lie algebra of type $C^{(1)}_\ell$ for $w=2\ell+1$. The conjecture is based on numerical evidence.
\end{abstract}
 \maketitle

{\it This note is dedicated to Marko Tadi\' c at the occasion of his 70th birthday.}
\smallskip

  \section{Introduction}

We write a partition of positive integer  $n$ in terms of frequencies  $f_j$---the number of occurrences of the part $j$ in the partition
\begin{equation}\label{E:partition}
n=\sum_{j\in\mathbb N} \, f_j\cdot j = \underset{\text{$f_1$
times}}{\underbrace{1+\dots+1}}+ \underset{\text{$f_2$ times}}{\underbrace{2+\dots+2}}+ \dots \, .
\end{equation}
It is clear that $f_j=0$ for all but finitely many $j\in\mathbb N$ and that the partition (\ref{E:partition}) is determined by its sequence of frequencies $(f_i\mid i\in\mathbb N)$.

The partition identities of Rogers (1894), Ramanujan (1913) and Schur (1917) for $k=1$, and the  partition identities of Gordon (1961) for $k\geq 2$, can be stated as:  
\smallskip

{\it Let $0\leq a\leq k$. The number of partitions of $n$ such that
\begin{equation}\label{E:RR difference conditions}
f_j+f_{j+1} \leq k\quad\text{for all $j$ and}
\end{equation}
\begin{equation}\label{E:RR initial conditions}
 f_1 \leq a\qquad\qquad
\end{equation}
equals the number of partitions of   $n$ into parts }
$\not\equiv 0, \pm (a+1)\mod{(2k+3)}$.
\smallskip

\noindent
The conditions (\ref{E:RR difference conditions}) on frequencies of two adjacent numbers are called {\it the difference conditions}, and the condition (\ref{E:RR initial conditions}) on the frequency of number $1$ is called {\it the initial condition}.
There are some other similar partition identities stating that the number of partitions of $n$ satisfying certain difference $\&$ initial conditions is equal to the number of partitions of $n$ with parts satisfying certain congruence conditions; these identities are often called the {\it classical} Rogers-Ramanujan type identities---see \cite{A}. 
On the other side, some parts of representation theory of affine Kac-Moody Lie algebras lead to Rogers-Ramanujan type {\it colored} partition identities. 

Let $N_1,\dots, N_r$, $r\geq2$, be non-empty subsets of the set of positive integers $\mathbb N$ and let $\mathcal N$ be the multiset
\begin{equation}\label{E:multiset calN}
\mathcal N=N_1\cup\dots\cup N_r.
\end{equation}
If a positive integer $a$ appears in several subsets $N_i$, then $a$ appears in the multiset $\mathcal N$ several times. To see these elements in $\mathcal N$ as different, for each positive integer $a$ 
we may ``color'' $a\in N_j$ with a ``color'' $j$ by writting $a_j=(a,j)\in N_j\times \{j\}$, and then write (\ref{E:multiset calN}) in terms of sets as
\begin{equation}\label{E:multiset calN as a set}
\mathcal N= (N_1\times \{1\})\cup\dots\cup  (N_r\times \{r\})\subset 
\mathbb N\times \{1,\dots, r\}.
\end{equation}
We say that elements in the multiset $\mathcal N$ appear in $r$ colors. In this note a colored partition of positive integer $n$ on the multiset $\mathcal N$ is
\begin{equation}\label{E:colored partition}
n=\sum_{a\in\mathcal N} \, f_a\cdot a .
\end{equation}
It is clear that $f_a=0$ for all but finitely many $a\in\mathcal N$ and that the partition (\ref{E:colored partition}) is determined by its ``sequence'' of frequencies $(f_a\mid a\in\mathcal N)$.
\begin{example}\label{Ex:1} Let $\mathcal N=N_1\cup N_2$, where
$$
N_1=\{j\in\mathbb N\mid j\equiv 2,8\mod{10}\},\quad
N_2=\{j\in\mathbb N\mid j\equiv 1,2,4,5,6,8,9\mod{10}\}.
$$
Then parts $a$ of colored partitions (\ref{E:colored partition}) for $\mathcal N=N_1\cup N_2$
appear in two colors, $1$ and $2$: parts $\equiv 2,8\mod{10}$ appear in both colors, and parts $\equiv 1,4,5,6,9\mod{10}$ appear only in  color $2$. Note that the generating function for colored partitions (\ref{E:colored partition}) is the infinite periodic product with modulus $10$:
\begin{equation}\label{E:infinite product}
\prod_{ j\equiv 1,2,2,4,5,6,8,8,9\mod{10}}(1-q^j)^{-1}.
\end{equation}
\end{example}

Lepowsky and Wilson gave in \cite{LW} a Lie theoretic interpretation of the classical Rogers-Ramanujan type partition identities in terms of characters of standard modules $L_{A^{(1)}_1}(\Lambda)$ for affine Kac-Moody Lie algebra of the type $A^{(1)}_1$. After their discovery it was expected that for each standard module $L_{{\mathfrak g(A)}}(\Lambda)$ for any affine Lie algebra ${\mathfrak g}(A)$ (cf. \cite{K}) there is a Rogers-Ramanujan type partition identity, where $A^{(1)}_1$ is just ``the smallest one'' on the list of all affine Lie algebras:
 \smallskip
    
$A^{(1)}_{1}, A^{(1)}_{\ell}, B^{(1)}_{\ell}, C^{(1)}_{\ell}, D^{(1)}_{\ell}, 
E^{(1)}_{6,7,8},  F^{(1)}_{4}, G^{(1)}_{2},
A^{(2)}_{2}, A^{(2)}_{2\ell}, A^{(2)}_{2\ell-1}, D^{(2)}_{\ell+1}, E^{(2)}_{6}, {D^{(3)}_{4}}. $
 \smallskip
 
 \noindent
 However, besides several sporadic results beyond $A^{(1)}_1$, so far this goal is not achieved.
 
 In \cite{CMPP} Rogers-Ramanujan type partition identities are conjectured for all standard $C^{(1)}_\ell$-modules, stating that the number of colored partitions of $n$ with parts satisfying certain congruence conditions is equal to the number of colored partitions (\ref{E:colored partition}) for a multiset $\mathcal N=\mathcal N_{2\ell+1}$ composed of $\ell$ copies of $\mathbb N$ and an additional copy of $(2\mathbb N+1)$,  satisfying difference $\&$ initial conditions similar to  (\ref{E:RR difference conditions})--(\ref{E:RR initial conditions}) , but much more complicated. Moreover, in  \cite{CMPP} another series of similar partition identities is conjectured for a multiset $\mathcal N=\mathcal N_{2\ell}$ composed of $\ell$ copies of $\mathbb N$, satisfying certain difference $\&$ initial conditions, but with no obvious connection to representation theory of affine Lie algebras.
 
In this note we conjecture yet another Rogers-Ramanujan type colored partition identities for a multiset $\mathcal N={\mathcal N}_{2\ell-1}^{\text{odd}}$, somewhat similar to the conjectured identities for standard $C^{(1)}_\ell$-modules, but again with no obvious connection to representation theory of affine Lie algebras.

 \section{Arrays with odd width $w$ and even first row}
 
 Let $\mathcal N={\mathcal N}_5^{\text{odd}}$  be the colored array of natural numbers with $5$ rows 
\begin{equation}\label{E:colored array of numbers}
\begin{matrix}
&{2_{1}}& &{4_{1}} & &{6_{1}} & &{8_{1}}  \\
{1_{1}}& &{3_{1}}& &{5_{1}}& &{7_{1}}&   \\
&{2_{2}} & &{4_{2}} & &{6_{2}} & &{8_{2}}   \\
{1_{2}}& &{3_{2}} & &{5_{2}} & &{7_{2}} &  \\
&{2_{3}} & &{4_{3}} & &{6_{3}} & &{8_{3}}   \\
\end{matrix}\quad\dots\, .
\end{equation}
$\mathcal N$ is a multiset composed of $2$ copies of $\mathbb N$ and an additional copy of $2\mathbb N$, but its elements are arranged in such a way that in  the first row are even numbers and that numbers increase by one going to the right on any diagonal. 
 
We consider colored partitions
\begin{equation}\label{E:colored partition on the array}
n=\sum_{a\in\mathcal N}f_a\cdot a,
\end{equation}
where $f_a$ is the frequency of the part $a\in\mathcal N$ in the colored partition (\ref{E:colored partition on the array}) of $n$. It is clear that  $f_a=0$ for all but finitely many $a\in\mathcal N$ and that the colored partition (\ref{E:colored partition on the array}) is determined by its array $\mathcal F$ of frequencies
\begin{equation}\label{E:array of frequencies}
\mathcal F= \quad
\begin{matrix}
& f_{2_{1}}& &f_{4_{1}} & &f_{6_{1}} & & f_{8_{1}}  \\
f_{1_{1}}& & f_{3_{1}}& & f_{5_{1}}& & f_{7_{1}}&   \\
&f_{2_{2}} & &f_{4_{2}} & &f_{6_{2}} & &f_{8_{2}} \\
f_{1_{2}}& &f_{3_{2}} & &f_{5_{2}} & &f_{7_{2}} &  \\
&f_{2_{3}} & &f_{4_{3}} & &f_{6_{3}} & &f_{8_{3}}   \\
\end{matrix} \quad\dots \, .
\end{equation}
We say that two elements in the array $\mathcal F$ are {\it adjacent} if they are simultaneously on two adjacent rows and two adjacent diagonals. For example, $f_{5_1}$ and $f_{7_1}$ in the second row are adjacent to $f_{6_1}$ in the first row and, just as well, adjacent to $f_{6_2}$ in the third row. We say that the set\footnote{
or the sequence
} $\{a_1,a_2,a_3,\dots\}$ is a downward path $\mathcal Z$ in the array $\mathcal F$ if $a_i$ is in the $i$-th row and if $(a_i, a_{i+1})$ is a pair of two adjacent elements for all $i$. For example, $\mathcal Z=\{f_{6_1},f_{5_1},f_{6_2},f_{7_2},f_{6_3}\}$ is a downward path in $\mathcal F$ and there are altogether $2^4$ downward paths through $f_{6_1}$ in the first row.

Let $k$ be a positive integer. We say that the frequency array $\mathcal F$ {\it satisfies level $k$ difference conditions} if
\begin{equation}\label{E:level k difference conditions}
\sum_{m\in\mathcal Z}m\leq k \quad\text{for all downward paths }\mathcal Z\text{ in } \mathcal F.
\end{equation}
Note that the level $k$ difference conditions for a frequency array $\mathcal F$ is similar to difference conditions (\ref{E:RR difference conditions}) for a sequence of frequencies $(f_i\mid i\in\mathbb N)$, but much more complicated.

Let $k_0, k_1, k_2, k_3\in\mathbb N_0$, $k=k_0+k_1+k_2+k_3>0$. We say that an array of frequencies $\mathcal F$ {\it is $(k_0, k_1, k_2, k_3)^{\text{odd}}$-admissible} if {\it the extended array of frequencies} 
\begin{equation}\label{E:extended array of frequencies}
\mathcal F\sp{(k_0, k_1, k_2, k_3)\text{odd}}=\quad
\begin{matrix}
&k_3&& f_{2_{1}}& &f_{4_{1}} & &f_{6_{1}} & & f_{8_{1}}  \\
k_2&&f_{1_{1}}& & f_{3_{1}}& & f_{5_{1}}& & f_{7_{1}}&   \\
&0&&f_{2_{2}} & &f_{4_{2}} & &f_{6_{2}} & &f_{8_{2}}   \\
k_1&&f_{1_{2}}& &f_{3_{2}} & &f_{5_{2}} & &f_{7_{2}} &  \\
&k_0&&f_{2_{3}} & &f_{4_{3}} & &f_{6_{3}} & &f_{8_{3}}   \\
\end{matrix}\quad\dots
\end{equation}
satisfies the level $k$ difference conditions, that is
\smallskip

\noindent
\begin{equation}\label{E:initial and difference conditions}
\sum_{m\in\mathcal Z}m\leq k \quad\text{for all downward paths }\mathcal Z\text{ in } \mathcal F^{(k_0, k_1, k_2, k_3)\text{odd}}.
\end{equation}
Note the difference between (\ref{E:level k difference conditions}) and (\ref{E:initial and difference conditions}): $(k_0, k_1, k_2, k_3)^{\text{odd}}$-admissible frequency array $\mathcal F$ satisfies the level $k$ difference conditions (\ref{E:level k difference conditions}), but in addition to that there are new conditions on the frequencies at the beginning of the array, somewhat similar to initial condition (\ref{E:RR initial conditions}), but much more complicated. For example, $f_{1_{1}}$ in the second row must be $\leq k_2$ because of (\ref{E:initial and difference conditions}) for the downward path $\mathcal Z=\{k_3, f_{1_{1}},0,k_1,k_0\}$.

We say that colored partitions (\ref{E:colored partition on the array}) with $(k_0, k_1, k_2, k_3)^{\text{odd}}$-admissible arrays of frequencies (\ref{E:array of frequencies}) are {\it $(k_0, k_1, k_2, k_3)^{\text{odd}}$-admissible colored partitions}.

Up till now we discussed the colored array of natural numbes  $\mathcal N={\mathcal N}_5^{\text{odd}}$ with $5$ rows, but all the notions can be extended to arrays   ${\mathcal N}_w^{\text{odd}}$ with $w=2\ell-1$ rows for $\ell=2, 3, 4, \dots$ with initial ``imaginary frequences'' being (from the bottom row to the top row):
\begin{itemize}
\item $(k_0, k_1, k_2)^{\text{odd}}=[k_0, k_1, k_2]$\qquad\quad for $\ell=2$, $w=3$,
\item $(k_0, k_1, k_2, k_3)^{\text{odd}}=[k_0, k_1, 0, k_2, k_3]$\qquad\quad for $\ell=3$, $w=5$,
\item $(k_0, k_1, k_2, k_3, k_4)^{\text{odd}}=[k_0, k_1, 0, k_2, 0, k_3, k_4]$\qquad\quad for $\ell=4$, $w=7$,
\item \qquad\qquad\qquad\dots \qquad
\item $(k_0, k_1,   \dots, k_{\ell-1}, k_\ell)^{\text{odd}}=[k_0, k_1, 0, k_2, 0, k_3, 0, \dots, 0, k_{\ell-1}, k_\ell]$ for $\ell\geq5$.
\end{itemize}
\smallskip

We say that colored partitions (\ref{E:colored partition on the array}) on $\mathcal N=\mathcal N_{2\ell-1}$
with  $(k_0, k_1,   \dots, k_{\ell-1}, k_\ell)^{\text{odd}}$-admissible arrays of frequencies (\ref{E:array of frequencies}) are 
{\it $(k_0, k_1,   \dots, k_{\ell-1}, k_\ell)^{\text{odd}}$-admissible colored partitions}.

\begin{conjecture}
Let $\ell\geq2$ and $k_0, k_1,\dots, k_\ell\in\mathbb N_0$, $k=k_0+\dots+k_\ell>0$.
Then the generating function for $(k_0, k_1,   \dots, k_{\ell-1}, k_\ell)^{\text{odd}}$-admissible colored partitions can be expresed as an infinite periodic product with modulus $2\ell+2k$.
\end{conjecture}
This conjecture is based on numerical evidence: we calculate\footnote{
by using a slightly modified code 21AAIC in \cite{CMPP} with built in option to choose even numbers in the top row (for p=0) or to choose odd numbers in the top row (for p=1);
} the number $a_n$ of $(k_0, k_1,   \dots, k_{\ell-1}, k_\ell)^{\text{odd}}$-admissible colored partitions of $n$ and then use Euler's factorization algorithm to write the generating function of partitions $\sum a_nq^n$ as an infinite periodic product---the Python code is available at \begin{verbatim}https://github.com/mirkoprimc/odd_w_odd
\end{verbatim}
Bellow are listed some results, where
\begin{itemize}
\item
``[1, 0, 0] product: [2, 3, 4] mod 6'' means that the conjectured generating function for
$(1, 0, 0)^{\text{odd}}$-admissible colored partitions is
$$
\prod_{ j\equiv 2,3,4\mod{6}}\frac{1}{(1-q^j)}
$$
\item
and ``[0, 1, 0] product: [1,2,-3,4,5] mod 6'' means that the conjectured generating function for
$(0, 1, 0)^{\text{odd}}$-admissible colored partitions is
$$
\frac{\prod_{ j\equiv 3\mod{6}}(1-q^j)}{\prod_{ j\equiv 1,2,4,5\mod{6}}(1-q^j)}.
$$
\end{itemize}
 \bigskip

\begin{verbatim}
[1, 0, 0] product: [2, 3, 4] mod 6

[0, 1, 0] product: [1, 2, -3, 4, 5] mod 6

[2, 0, 0] product: [2, 3, 4, 4, 5, 6] mod 8

[1, 1, 0] product: [1, 2, 4, 4, 6, 7] mod 8

[1, 0, 1] product: [2, 2, 3, 5, 6, 6] mod 8

[0, 2, 0] product: [1, 2, 2, 6, 6, 7] mod 8

[3, 0, 0] product: [2, 3, 4, 4, 5, 6, 6, 7, 8] mod 10

[2, 1, 0] product: [1, 2, 4, 4, 5, 6, 6, 8, 9] mod 10

[2, 0, 1] product: [2, 2, 3, 4, 5, 6, 7, 8, 8] mod 10

[1, 2, 0] product: [1, 2, 2, 4, 5, 6, 8, 8, 9] mod 10

[1, 1, 1] product: [1, 2, 3, 4, 4, -5, 6, 6, 7, 8, 9] mod 10

[0, 3, 0] product: [1, 2, 2, 3, 4, -5, 6, 7, 8, 8, 9] mod 10

[4, 0, 0] product: [2, 3, 4, 4, 5, 6, 6, 7, 8, 8, 9, 10] mod 12

[3, 1, 0] product: [1, 2, 4, 4, 5, 6, 6, 7, 8, 8, 10, 11] mod 12

[3, 0, 1] product: [2, 2, 3, 4, 5, 6, 6, 7, 8, 9, 10, 10] mod 12

[2, 2, 0] product: [1, 2, 2, 4, 5, 6, 6, 7, 8, 10, 10, 11] mod 12

[2, 1, 1] product: [1, 2, 3, 4, 4, 6, 6, 8, 8, 9, 10, 11] mod 12

[2, 0, 2] product: [2, 2, 3, 4, 4, 5, 7, 8, 8, 9, 10, 10] mod 12

[1, 3, 0] product: [1, 2, 2, 3, 4, 6, 6, 8, 9, 10, 10, 11] mod 12

[1, 2, 1] product: [1, 2, 2, 4, 4, 5, 7, 8, 8, 10, 10, 11] mod 12

[0, 4, 0] product: [1, 2, 2, 3, 4, 4, 8, 8, 9, 10, 10, 11] mod 12

[1, 0, 0, 0, 0] product: [2, 3, 4, 5, 6] mod 8

[0, 1, 0, 0, 0] product: [1, 2, 4, 6, 7] mod 8

[2, 0, 0, 0, 0] product: [2, 3, 4, 4, 5, 5, 6, 6, 7, 8] mod 10

[1, 1, 0, 0, 0] product: [1, 2, 3, 4, 4, 6, 6, 7, 8, 9] mod 10

[1, 0, 0, 1, 0] product: [1, 2, 2, 4, 5, 5, 6, 8, 8, 9] mod 10

[1, 0, 0, 0, 1] product: [2, 2, 3, 3, 4, 6, 7, 7, 8, 8] mod 10

[0, 2, 0, 0, 0] product: [1, 2, 2, 3, 4, 6, 7, 8, 8, 9] mod 10

[0, 1, 0, 1, 0] product: [1, 1, 2, 4, 4, 6, 6, 8, 9, 9] mod 10

\end{verbatim}
\begin{remark} From the list above we see that we may expect Roger-Ramanujan type colored partition identities for most of $(k_0, k_1,   \dots, k_{\ell-1}, k_\ell)^{\text{odd}}$-admissible para\-me\-ters---like the conjectured product formula (\ref{E:infinite product}) for the generating function of $(1, 2, 0)^{\text{odd}}$-admissible colored partitions.
It seems that for all the other parameters there is no infinite periodic product; the first possible case are parameters $[0, 0, 1, 0, 0]$ for $\ell=3$ for which our code gives
\begin{verbatim}
/oddWodd/allWcasesProd
the first row parity p = 0 , the highest_weight = [0, 0, 1, 0, 0] , 
N = 18
the exponents of the conjectured periodic product: 
[2, 3, 3, 4, 4, -6, -6, -7, -7, 8, 9, 9, 9, 9, 10, 10, 10, 10, 10, 
10, 10, 11, 11, -12, -12, -12, -12, -12, -12, -12, -13, -13, 
-13, -13, -13, -13, -13, -13, -13, -13, -13, -13, -14, -14, 
-14, -14, 15, 15, 15, 15, 15, 15, 15, 15, 15, 15, 15, 15, 15, 
15, 16, 16, 16, 16, 16, 16, 16, 16, 16, 16, 16, 16, 16, 16, 
16, 16, 16, 16, 16, 16, 16, 16, 16, 16, 16, 16, 16, 16, 16, 
17, 17, 17, 17, 17, 17, 17, 17, 17, 17, 17, 17, 17, 17, 17, 
17, 17, 17, 17, 17, -18, -18, -18, -18, -18, -18, -18, -18, 
-18, -18, -18, -18, -18, -18, -18, -18, -18, -18, -18, -18, 
-18, -18, -18]
\end{verbatim}
The list above encodes Euler's product for the first $19$ terms in the generating function of $[0, 0, 1, 0, 0]$-admissible partitions (in the sense of \cite{CMPP}) for $p=0$:
$$
\frac{(1-q^6)^2(1-q^7)^2(1-q^{12})^7(1-q^{13})^{12}(1-q^{13})^{12}(1-q^{14})^{4}\dots}
{(1-q^2)(1-q^3)^2(1-q^4)^2(1-q^8)(1-q^9)^4(1-q^{10})^7(1-q^{11})^2(1-q^{15})^{14}\dots}.
$$
\end{remark}

\section*{Acknowledgement}
I thank Marko Tadi\' c for us being friends and the members of our Representation theory seminar in Zagreb for so many years.

I thank Stefano Capparelli, Arne Meurman and Andrej Primc for many of their ideas that are interwoven in this note.

This work is partially supported by Croatian Science Foundation under the project 8488 and by the QuantiXLie Centre of Excellence, a project cofinanced by the Croatian Government and European Union through the European Regional Development Fund - the Competitiveness and Cohesion Operational Programme (Grant KK.01.1.1.01.0004).



\begin{thebibliography}{CMPP}

\bibitem[A]{A}
G. E. Andrews,
\textit{ The theory of partitions},
Encyclopedia of Mathematics and Its Applications, Vol. 2, Addison-Wesley, 1976.

\bibitem [CMPP]{CMPP}   S. Capparelli, A. Meurman, A. Primc and M. Primc
\textit{New partition identities from $ C^{(1)}_{\ell}$-modules}, Glasnik Matemati\v cki  \textbf{57} (2022),161--184.
335--355.

\bibitem [K]{K}
V. G. Kac, \textit{Infinite-dimensional Lie algebras} 3rd ed,
Cambridge Univ. Press, Cambridge, 1990.

\bibitem [LW]{LW}
J. Lepowsky and R. L. Wilson,  \textit{The structure of standard
modules, I: Universal algebras and the Rogers-Ramanujan
identities}, Invent. Math. \textbf{77} (1984),  199--290;
\textit{II: The case $A_1^{(1)}$, principal gradation}, Invent.
Math. \textbf{79} (1985),  417--442.




\end{thebibliography}
\end{document}